\documentclass{amsart}
\usepackage{graphicx}
\usepackage{amsfonts}
\usepackage{amssymb}
\usepackage{amsmath}
\usepackage{t1enc, epsfig}
\usepackage[latin1]{inputenc}
\usepackage[english]{babel}
\usepackage{pstricks}
\usepackage{pst-plot}
\usepackage{multido}

\newcommand{\RR}{\mathbb{R}}
\newcommand{\ZZ}{\mathbb{Z}}

\newcommand{\Rz}{\mathbb{R}^{\mathbb{Z}}}
\newcommand{\co}[1]{o(#1)}

\newcommand{\ww}[1]{\widetilde{#1}}
\newcommand{\cI}{ {\mathcal I}}

\newtheorem{thedec}{The Decomposition Theorem}[section]
\newtheorem{definition}{Definition}[section]
\newtheorem{lemma}{Lemma}[section]

\begin{document}

%\begin{frontmatter}

\title[Closed and exact functions]{Closed and exact functions in the context of  Ginzburg-Landau models}
\author{BY Anamaria Savu}
%\author[Savu]{BY Anamaria Savu} \corauthref{cor} \thanksref{tha}}
%\thanks[tha]{Research supported by the Ontario Graduate Studies Fellowship and the University of Toronto Fellowship}
%\corauth[cor]{Corresponding author.}
%\address[Savu]{Department of Mathematics and Statistics, Queen's University, Kingston, ON, K7L 3N6, Canada}
%\ead{ana@mast.queensu.ca}

\address{Department of Mathematics and Statistics, Queen's University, Kingston, ON, K7L 3N6, Canada}
\email{ana@mast.queensu.ca}

%\subjclass{Primary 42B05; secondary 81Q50, 42A16.}
%\date{ 27th Feb. 2005}
%\commby{XXXX}

\begin{abstract}
For a general vector field we exhibit two Hilbert spaces, namely 
the space of so called ``closed functions'' and the space of ``exact functions'' 
and we calculate the codimension of the space of ``exact functions'' 
inside the larger space of ``closed functions''.  
In particular we provide a new approach for the known cases: 
the Glauber field and the second-order Ginzburg-Landau field
and for the case of the fourth-order Ginzburg-Landau field.
% that are relevant for  statistical 
%physics questions like: making rigorous sense of the 
%fluctuation-dissipation equation and finding the large scale dynamics
% of the interacting particle systems, the Glauber model, 
%the second-order Ginzburg-Landau model 
%and the continuum solid-on-solid model, also known as the fourth-order 
%Ginzburg-Landau model. All models are assumed to have quadratic potential.  
%The result in the case of Glauber field and second-order Ginzburg-Landau
%field were previously known. We present a novel approach, based on 
%Fourier analysis, that allows us to calculate the codimension
%of the space of ``exact functions'' inside the space of ``closed functions''
%for any vector field. 

%The result for the continuum solid-on-solid  model is new and 
%can be used to find the large scale dynamics of this system, see \cite{Sav}. 
%A nongradient model, like continuum solid-on-solid model has 
%the Hilbert space of closed functions associated with. A 
%natural subspace of the space of closed functions is the 
%space of exact functions. In the present paper we define the
%space of exact functions and the space of closed functions
%for the continuum solid-on-solid model and we show that the 
%space of exact functions has codimension one inside the 
%space of closed functions.
\end{abstract}

%\begin{keyword}
%\end{keyword}

%\end{frontmatter}

\maketitle

\tableofcontents
\newpage
%%%%%%%%%%%%%%%%%%%%%%%%%%%%%%%%%%%%%%%%%%%%%%%%%%%%%%%%%%%
%\section{Introduction}
%%%%%%%%%%%%%%%%%%%%%%%%%%%%%%%%%%%%%%%%%%%%%%%%%%%%%%%%%%%
\section{Introduction} \label{intro}

Statistical physics has developed a whole variety of interacting 
particle systems that capture some aspects of the movement of 
particles on the microscale. An interacting particle system
is usually a complex Markov process with a finite or infinite
state space. By taking an appropriate scaling limit of an
interacting particle system we expect to derive the evolution
of the system on the macroscale, in general a nonlinear 
partial differential equation. 
It is fairly well understood the transition from the microscopic 
scale to macroscopic scale, at least for some systems, and 
in this notes we take for granted this step.

%There has been  a great deal of 
% investigation (see \cite{Kip}, or \cite{Spo}),
%to explain how the transition from the microscopic scale to macroscopic 
%scale is achieved.   
  
The most interesting microscopic models constructed so far, lack
the so called gradient condition. This condition corresponds to 
the Fick's law of fluid dynamics according to which the instantaneous 
current $w$ of particles over a bond is the gradient
$\tau h - h$ of some local function.
Since the work of Varadhan \cite{Var1}, Quastel \cite{Qua1} 
and Varadhan and Yau \cite{Var2}
on nongradient systems, new ideas have been imposed in the field. 
The main idea is that a nongradient system has a generalized 
Fick's law, also called the fluctuation-dissipation equation, 
of the form
$$w \approx \hat{a}(m)(\tau h -h) + L g, $$
where $\hat{a}$ is the transport coefficient depending on the 
particle density $m$ in a microscopic cube, $h$ is some local 
function, and $L$ is the generator of the
microscopic dynamics. The $Lg$ part of the approximate
equation above is negligible on the macroscale, and is called
the fluctuation part of the equation. 

One of the main difficulty in finding the scaling limit of a 
nongradient system is to make rigorous sense of the 
fluctuation-dissipation equation. As it has been shown in 
\cite{Var1}, \cite{Qua1}, 
\cite{Var2} the current w, the gradient
$\tau h -h$ and the fluctuations $Lg$ are elements of the Hilbert 
space of ``closed functions'' and the fluctuation-dissipation equation 
is a consequence of a direct-sum decomposition of this Hilbert space. 
The gradient part 
$\tau h-h$ of the current $w$ that survives after taking the 
scaling limit of the model is just the projection of $w$ onto
a one dimensional subspace of the Hilbert space of ``closed functions''. 
The remaining negligible fluctuations $Lg$ are vectors of the Hilbert 
subspace of ``exact function''. 

The purpose of the present paper is not to show how the Hilbert 
space of ``closed functions'' and ``exact functions'' arises in 
the context of interacting particle systems, but rather to motivate 
the direct-sum decomposition of the Hilbert space of ``closed functions''
and to find the codimension of the space of ``exact functions'' inside
the space of ``closed functions''. We calculate this codimension 
for an arbitrarily chosen vector field. 
The three continuum models known as the Glauber system, the 
second-order Ginzburg-Landau system and the continuum solid-on-solid 
model, also called the fourth-order Ginzburg-Landau system are 
covered by our general result. 
Our approach to establish the direct-sum decomposition of the 
Hilbert space of ``closed functions'' is new and differs from
the approach used before to study the first two models 
(see Varadhan \cite{Var1}). We have followed a different path 
based on Fourier analysis that has allowed us to handle a 
general vector field. 

%Since we  have tried unsuccessfully 
%to extend Varadhan's method \cite{Var1}, we follow a different path 
%based on Fourier analysis that allow us to handle a general vector field. 

%%%%%%%%%%%%%%%%%%%%%%%%%%%%%%%%%%%%%%%%%%%%%%%%%%%%%%
%The decomposition theorem
%%%%%%%%%%%%%%%%%%%%%%%%%%%%%%%%%%%%%%%%%%%%%%%%%%%%%%
\section{The decomposition theorem} \label{decomp}

In this section we introduce some terminology and state the 
main result.

The Hermite polynomials provide an orthogonal basis for the Hilbert space
of functions defined on the real axis, that are square integrable
with respect to the Gaussian probability measure 
$\frac{1}{\sqrt{2\pi}}\exp(-\frac{x^2}{2}) dx$.   
The $i$th Hermite polynomial is defined through
$$ H_i(x) =  \frac{(-1)^i}{i!} \exp\bigg(\frac{x^2}{2}\bigg) 
\bigg(\frac{d^i}{dx^i}\exp\bigg(-\frac{x^2}{2}\bigg) \bigg), 
\quad i \in \mathbb{N}.$$
We stress that $H_i$ is not normalized to have $L^2$ norm 1 
with respect to the probability measure 
$\frac{1}{\sqrt{2\pi}}\exp(-\frac{x^2}{2}))dx$, but rather 
$\frac{1}{\sqrt{i!}}.$ 

There is an extension of Hermite polynomials to more variables. 
A multi-index is a double-sided infinite sequence $I = \{ i_n \}_{n \in \ZZ}$
of positive integers, with at most finitely many non-zero entries. The degree 
of a multi-index is $|I|= \sum_{n \in \ZZ} i_n$. 
Call $\mathcal{I}$ the set of multi-indices and $\mathcal{I}_N$ the set of 
multi-indices of fixed degree $N$. 
The multidimensional Hermite polynomials are
$$ H _I(x) = \Pi_{n\in \mathbb{Z}} H_{i_n}(x_n), \quad I \in \cI.$$ 
We make the convention that if a multi-index $I$ has some strictly 
negative entries then $H_I =0$. Together the multidimensional Hermite 
polynomials, $\{ H_I \}_{I \in \mathcal{I}}$ form an orthogonal basis 
for the Hilbert space of functions defined on $\Rz$, that are square 
integrable with respect to the probability measure 
$$d\nu_{0}^{gc} =\bigotimes_{i \in \ZZ} \frac{1}{\sqrt{2\pi}}\exp\bigg(-\frac{x_i^2}{2}\bigg)dx_i.$$
It is interesting to note that this Hilbert space is a model for the symmetric Fock space over
the space of square summable, double-sided sequences $l^2(\ZZ)$, 
%$\Gamma_s(l^2(\ZZ))$, 
and decomposes as a direct sum of the degree
$N$ subspaces 
$${\mathcal H}_N = \{ H_I \; | \; |I|=N \}^c.$$
%\simeq l^2(\ZZ)^{\odot N}$.
%$$L^2(\Rz, d\mu ) = \bigoplus_{N \geq 0} {\mathcal H}_N, $$
The superscript on the line above, means that we take the closed linear span of the set.

The shift $\tau$ acts on configurations as $(\tau (x ))_{n} =x_{n+1}$ and 
on functions as $(\tau f)(x) = f(\tau x)$.  $\tau^n$ stands for the $n$-fold 
composition $\tau \circ \dots \circ \tau$.
If a multi-index $I=(i_n)_{n\in \mathbb{Z}}$ has $i_n=0$ for all $n <0$ we shall say 
that the multi-index is supported on the set of positive integers. 
We shall use the notation $\delta_n$ for the multi-index that corresponds to 
the configuration with a single particle at the site $n$. 
Two multi-indices can be added and the addition is point-wise. 

The action of the annihilation, creation, and shift operators on the multidimensional
Hermite polynomial $H_I$ is very simple:
$$ \partial_{n} H_I(x) = H_{I-\delta_n} (x), \quad
   (x_n - \partial_{n} ) H_I(x) = H_{I+\delta_n} (x) , \quad
   \tau H_I = H_ {\tau^{-1}I } .$$
Above $\partial_n$ stands for the partial derivative with respect 
to the $n$th coordinate. 

Given a double-sided sequence of real numbers $(a_k)_{k \in \ZZ}$, that
are all but finitely many zero  we introduce the vector field 
$D_0=\sum_{k \in \ZZ} a_k \partial_k$ with constant coefficients. 
Translating $a$'s to the left or to the right produces a new 
sequence that defines the vector field $D_n=\sum_{k \in \ZZ} a_k \partial_{k+n}$,
$n \in \ZZ$. Now we have the setup needed to introduce the closed 
and exact functions.
%{\bf The space of closed and exact functions.}
\begin{definition}
We shall say that a function $\xi \in L^2(\Rz, d\nu_{0}^{gc})$ is 
closed (or more precisely, $D_0$-closed) if it satisfies in the weak sense
\begin{equation} \label{closedt}
D_n (\tau^m \xi) = D_m (\tau^n \xi)
\end{equation}
for all integers $m$ and $n$.
Let ${\mathcal C}_D$ denote the space of all $D_0$-closed functions.
\end{definition}
\begin{definition}
We shall say that a function $\xi^g \in L^2(\Rz, d\nu_{0}^{gc})$ is
exact (or more precisely $D_0$-exact) if there is a local function 
$g$, a function that depends on finitely many co-ordinates, such that
\begin{equation} \label{closed}
\xi^g = D_0 \bigg( \sum_{k \in \mathbb{Z}} \tau^k g \bigg) =
\sum_{k \in \mathbb{Z}} D_0(\tau^k g).
\end{equation}
Let ${\mathcal E}_D$ denote the closed linear span of the set of $D_0$-exact 
functions.
\end{definition}
Although the infinite sum $\sum_{k \in \mathbb{Z}} \tau^k g$  does not 
make sense, after applying the differential operator $D_0$ we get a 
meaningful expression. Since $g$ is a local function, the vector field 
$D_0$ kills all but finitely many terms of the infinite formal sum. 

The terminologies of exact and closed functions are not arbitrarily chosen. 
We can define formally the form $w = \sum_{n \in \ZZ} \tau^n \xi \; dx_n$
and the boundary operator $ df = \sum_{n \in \ZZ} D_n (f)\; dx_n.$ 
It is not hard to see, with these new definitions, that the form $w$ 
is closed ($dw=0$), in the vector calculus sense, if and only if 
$D_n (\tau^m \xi) = D_m (\tau^n \xi)$, i.e., if and only if $\xi$ 
is a closed function.

Knowing that any exact function is closed a natural question to ask is about 
the codimension of the space of exact functions inside the space of closed 
functions.  In this paper we provide the answer for this question.

\begin{thedec}  \label{the1}
Let $D_0=\sum_{k \in \ZZ}a_k \partial_k$ be a vector field with constant
real coefficients. All but finitely many numbers in the sequence 
$(a_k)_{k \in \ZZ}$ are zero. The following decomposition results hold:
\begin{itemize}
\item[a)] If the sum of the coefficients of the vector field $D_0$ is not
equal to zero then $${\mathcal C}_D =  {\mathcal E}_D.$$
\item[b)] If the sum of the coefficients of the vector field $D_0$ is 
equal to zero then  $${\mathcal C}_D = \mathbb{R}\bf{1}\oplus {\mathcal E}_D.$$
\end{itemize}
\end{thedec}

{\bf Idea of the proof for the decomposition theorem \ref{the1}.} We outline the 
main ideas used to prove the decomposition theorem. 
We shall show later that a function $\xi$ is $D_0$-closed if and only if 
the projections $\mathrm{Proj}_{\mathcal{H}_N}\xi$, $N \geq 0$ are $D_0$-closed.
Degree $0$ subspace is easy to analyze since it is one dimensional. Any constant 
function is always $D_0$-closed, but is exact if and only if the sum of the 
coefficients of $D_0$ is not equal to zero. If the sum of the coefficients of the 
$D_0$ is equal to zero, then any $D_0$-closed function is orthogonal on the 
degree $0$ subspace. Therefore the result of the theorem holds if we can prove 
that a given $D_0$-closed function $\xi$ in $\mathcal{H}_N$, $N \geq 1$, the 
function $\xi$ can be approximated with $D_0$-exact functions. 
 
We shall investigate the properties of the Fourier coefficients of closed and exact 
functions, and we shall rather establish that the Fourier coefficients of a closed 
functions can be approximated in the appropriate sense with Fourier coefficients 
of exact functions. The ideas will be elaborated in the following sections.   

{\it Note.} In two cases relevant for statistical physics questions, 
namely the second-order Ginzburg-Landau vector field $Y_0=\partial_1-\partial_0$ 
and the fourth-order Ginzburg-Landau vector field $X_0=\partial_1-2\partial_0+\partial_{-1}$, 
the decomposition result of Theorem \ref{the1} is equivalent with the fluctuation-dissipation 
equation mentioned in the introduction section of the paper. 

{\it Note.} To get a flavor of the result stated in Theorem \ref{the1} 
we give some examples of exact and closed functions in the case of the fourth-order 
Ginzburg-Landau field, $X_0=\partial_1-2\partial_0+\partial_{-1}$:  $x_n + x_{-n}-2x_0$ are $X_0$-exact,
 ${\bf 1}$, $x_0$, $x_n+x_{-n}$ are examples of $X_0$-closed 
but not $X_0$-exact functions. A strange phenomena appears for: besides the function 
{\bf{1}} there exists another function that is $X_0$-closed and not $X_0$-exact, namely 
$x_0$. Therefore one may expect that the codimension of the space of exact 
functions is two. This is not the case and $x_0$ can be approximated with 
exact functions.

%%%%%%%%%%%%%%%%%%%%%%%%%%%%%%%%%%%
%\section{The set of multi-indices}
%%%%%%%%%%%%%%%%%%%%%%%%%%%%%%%%%%%
\section{The set of multi-indices} \label{multiind}

A multi-index $I = \{i_n\}_{n \in \ZZ}$ can be thought of as 
a configuration of particles sitting on the sites of the lattice $\ZZ$. On top of the site $n$ sit
$i_n$ particles. Rather than saying how many particles are at each site, we give the positions
of the particles. This way we obtain a vector 
\begin{equation} \label{coding}
z_I = ( \underbrace{ n_1 \dots n_1}_{i_{n_1}}, \dots ,  \underbrace{ n_k \dots n_k}_{i_{n_k}}).
\end{equation} 
that lists, in increasing order, all occupied sites of $I$ repeated according to the number of 
particles that occupy the site. We assume the only non-zero entries of the multi-index $I$
are $i_{n_1}, \dots i_{n_k}$. Note that the dimension of the vector $z_I$ is the degree of 
the multi-index $I$. If the multi-index has zero degree then $z_I$ is just a point. We say 
that $z_I$ is a new coding of the multi-index $I$. This correspondence shows that the
set $\cI_N$ is bijective with the set of vectors of $\ZZ^N$ with entries in increasing 
order or is in bijection with the quotient space $\ZZ^N/S_N$, where $S_N$ is the group of 
permutations of $N$ letters. 

For the results that follow we need to say more about the set of multi-indices. 
We partition the set of multi-indices into orbits with the help of the 
group action 
\begin{equation} \label{action}
 \ZZ \times \ZZ^{\ZZ} \longrightarrow \ZZ^{\ZZ} \quad 
(n, I) \longmapsto n \cdot I :=  \tau^n(I -\delta_n + \delta_0).
\end{equation}
When restricted to $\ZZ \times \cI$ the map (\ref{action}) is not an action any 
more since the multi-indices that enumerate the basis of the $L^2$ space are constrained to have  
positive entries. 

The orbits of the action (\ref{action}) provide a partition of the set of 
multi-indices $\ZZ^\ZZ$. For each multi-index $I \in \mathcal{I} $ we define $\co{I}$ to be 
the shadow of the orbit of $I$ on the set $\cI$, i.e., 
$\co{I} = \{ J | J= n \cdot I  \quad n \in \ZZ \} \cap \cI = \{ J | J= n \cdot I\quad n\in s(I)\}.$
Here, $s(I)= \{ n \in \ZZ \; | \; i_n \neq 0 \}$ is the finite set of occupied positions of $I$. 
From now on we will refer to $\co{I}$ as the orbit of $I$, although this is just a part of 
the actual orbit of the action. It has the advantage of being finite since the multi-index 
$I$ has all but finitely many entries zero and there are just finitely many $n$'s that after 
acting on $I$ give rise to a multi-index with positive entries. All the multi-indices in 
the same orbit have the same degree. 
The orbits partition $\cI$  and $\cI_N$. Call $\mathcal{O}$ the set of orbits, and 
$\mathcal{O}_N$ the set of orbits containing multi-indices of degree $N$.

It is worth mentioning that inside each orbit $\co{I}$ there exists a unique representative 
supported on the positive integers. Denote this multi-index by $R(\co{I})$. To see that this 
is true let us assume that $I$ has some particles in some negative position, i.e., 
there exists some $n < 0$ such that $i_n \geq 1$. If $k$ is the 
leftmost occupied position of $I$ and $k<0$, then $k \cdot I \in \co{I}$ is supported on 
the positive integers. Let us call $\mathcal{R}$ the set of all representatives, and 
$\mathcal{R}_N$ the set of degree-$N$ representatives.

%{\bf The orbit space.} 
So far the orbit space $\{\co{I}\}_{I\in \cI}$ is an abstract object. 
Fortunately we are able to give a concrete description of the orbit space.  
For this purpose it is very useful to know that each orbit, $o$ has a unique representative 
$R(o)$ supported on the positive semi-axis. The vector $z_{R(o)}$ is a point in the positive 
cone $\mathcal{C}_N^+=\{ z \in \ZZ^N | z=(z_1, \dots, z_N), \; 0 \leq z_1 \leq \dots \leq z_N \}$.
Therefore the set of representatives, and in particular the set of orbits $\mathcal{O}_N$, 
are bijective with the cone $\mathcal{C}_N^+$. Since there is only one multi-index 
with zero degree,  ${\bf 0}=(0)_{n \in \ZZ}$, the sets $\cI_0$, ${\mathcal O}_0$ and 
${\mathcal R}_0$ contain just a single element. By convention, ${\mathcal C}_0^+$ is just
 one-point set. 

We can say even more about this picture. The cone $\mathcal{C}_N^{+}$, itself is an orbit space, 
which we shall describe below.

Let us define the transformations that acts on the lattice $\ZZ^N$, for any $1\leq i,j \leq N$,
\begin{equation}
\sigma_{i,j}:\ZZ^N \rightarrow \ZZ^N, \quad 
\sigma_{i,j}(z_1,\dots ,z_i, \dots, z_j \dots ,z_N) = (z_1,\dots ,z_j, \dots, z_i, \dots ,z_N)
\end{equation}
$$ \mathrm{and} \quad \gamma_1: \ZZ^N \rightarrow \ZZ^N, \quad
\gamma_1 (z_1,z_2,\dots ,z_N) = (-z_1, z_2-z_1\dots ,z_N-z_1). $$
The smallest group generated by $\sigma_{i,j}$, $1 \leq i,j \leq N$ and $\gamma_1$ will be 
denoted by $\widetilde{S}_N$. To see that $\ww{S}_N$ is isomorphic with the group of 
permutations of $N$ letters, we write down the basic relations among the generating transformations: 
$(\gamma_1 \sigma_{1,2})^3 = \bf{id}$ and $\gamma_1 \sigma_{i,i+1}=\sigma_{i+1,i}\gamma_1$, 
$1 \leq i \leq N-1$.
The group $\ww{S}_N$ has $S_N$, the group of permutations of $N$ letters as subgroup, and
 $\ww{S}_N$ decomposes into left cosets with respect to $S_N$, as $\ww{S}_N = S_N \cup 
\gamma_1 S_N \dots \cup \gamma_N S_N $, where the transformations $\gamma_i$ are 
\begin{equation}
\gamma_i: \ZZ^N \rightarrow \ZZ^N, \quad \gamma_i (z_1,z_2,\dots ,z_N) = (-z_i, z_2-z_i\dots ,z_N-z_i). 
\end{equation}

It is interesting to note that $\ZZ^N/\ww{S}_N$ is bijective with the cone $\mathcal{C}_N^+$, 
as the next argument proves. 
Any orbit of $\ZZ^N/\ww{S}_N$ contains at least one vector, let say $z$, with components
in increasing order. If this vector does not have positive co-ordinates, it means $z_1 < 0$.
 But $(-z_1, z_2-z_1, \dots , z_N-z_1)$ is still a point in the orbit of $z$ under the 
action of $\ww{S}_N$. We can rearrange the coordinates of the new vector to be in increasing 
order and hence the orbit of $z$ under the action of $\ww{S}_N$ contains at least one vector 
of the cone $\mathcal{C}_N^+$. To see that the orbit of $z$ does not contain more than 
one vector of $\mathcal{C}_N^+$ we use the coset decomposition of $\ww{S}_N$. If $z$ 
is in $\mathcal{C}_N^+$, then rearranging the co-ordinates of $z$ we obtain either
the vector $z$ or some vector outside the cone $\mathcal{C}_N^+$. If we act on $z$ 
or some other vector obtained from $z$ by changing the places of the co-ordinates, 
with either of the transformations $\gamma_1, \dots \gamma_N$ we get a vector that
has at least one negative co-ordinate, so does not belong to $\mathcal{C}_N^+$. 

Now we can say that the set of orbits $\mathcal{O}_N$ is bijective with the cone 
$\mathcal{C}_N^+$, and hence with the quotient space $\ZZ^N/\ww{S}_N$. The bijection is
$ o \in \mathcal{O}_N \mapsto z_{R(o)} \in \mathcal{C}_N^+$.

In addition, if $I$ and $J$ are two multi-indices in the same orbit of the action 
(\ref{action}) then $z_I$ and $z_J$ are in the same orbit of the action of $\ww{S}_N$ on $\ZZ^N$.
Assume that $J=n_j \cdot I$ with $I= \sum_{i=1}^k a_{i} \delta_{n_i}$, where $a_{i} \neq 0$ 
and $n_1 \leq \dots \leq n_k$. Then $J = \sum_{i=1,...,k, i\neq j} a_{i} \delta_{n_i -n_j}+
(a_{j}-1)\delta_0 + \delta_{-n_j}$, and so
$$ z_I = ( \underbrace{ n_1, \dots ,n_1}_{a_1}, \dots ,  \underbrace{ n_k, \dots ,n_k}_{a_k})$$
$$ z_J = ( -n_j,\underbrace{ n_1-n_j, \dots ,n_1-n_j}_{a_1}, \dots , \underbrace{ 0, \dots ,0}_{a_j-1}, 
\dots, \underbrace{ n_k -n_j, \dots ,n_k-n_j}_{a_k}) .$$
It follows that $z_J$ is the image of $z_I$ under some element of $\ww{S}_N$. 

We shall denote by $z \stackrel{S_N}{\sim} z'$ and $z \stackrel{\ww{S}_N}{\sim} z'$
two lattice points $z$ and $z'$ that have the same image in the quotient space 
$\ZZ^N/S_N$ and $\ZZ^N/\ww{S}_N$, respectively.

Before we leave this section it is important to notice the following crucial facts.
Let $N \geq 1$. Since $\mathcal{I}_N$ is identified with $\ZZ^N/S_N$ we can think of any 
function $\hat{\xi}:\mathcal{I}_N \to \RR$ as being a $S_N$-invariant function
$\hat{\xi}:\ZZ^N \to \RR$, where $ \hat{\xi}(z)=\hat{\xi}(I)$ if $z \stackrel{S_N}{\sim} z_I$.
Similarly, since $\mathcal{O}_N$ is identified with $\ZZ^N/\ww{S}_N$ we can 
think of any function $c:\mathcal{O}_N \to \RR$ as being a $\ww{S}_N$-invariant
function $\ww{c}:\ZZ^N \to \RR$, where $\ww{c}(z) = c(o)$ if there exists a multi-index 
$I \in o$ such that $z \stackrel{S_N}{\sim} z_I$.

%%%%%%%%%%%%%%%%%%%%%%%%%%%%%%%%%%%%%%%%
%\section{Properties of closed functions and of exact functions}
%%%%%%%%%%%%%%%%%%%%%%%%%%%%%%%%%%%%%%%%
\section{Properties of closed functions and of exact functions} \label{clos-exac}

This section contains a detailed study of closed and exact functions. 

\textbf{Closed functions.} We start with 
a very simple but important property of closed functions.
\begin{lemma} \label{lema1}
Assume $D_0$ is a vector field with constant coefficients, $D_0=\sum_{k \in \ZZ} a_k \partial_k$.
 All but finitely many coefficients of the vector field $D_0$ are zero. 
A function $\xi \in L^2(\RR^{\ZZ}, d\nu_0^{gc})$ is $D_0$-closed if and only if the 
projection $Proj_{{\mathcal H}_N} \xi$ onto the degree $N$  subspace $\mathcal{H}_N$ is 
 $D_0$-closed, for any $N \geq 0$.
\end{lemma} 
\noindent
{\bf Proof.}
We denote by $\partial_j$ the differential operator with respect to the 
$j^{\mathrm{th}}$ coordinate, and by $\partial_j^*= -\partial_j +x_j$ the 
adjoint operator of $\partial_j$. The adjoint is taken with respect to the 
inner product $<,>$ of $L^2(\RR^{\ZZ}, d\nu_0^{gc})$. 
The operators $\partial_j$ and $\partial_j^*$ are bounded operators when 
restricted to a degree subspace, although they are unbounded on the 
whole $L^2(\RR^{\ZZ}, d\nu_0^{gc})$ space. 

If $\xi \in \mathcal{H}_N$, with Fourier series $\xi = \sum_{I \in \mathcal{I}_N} 
\hat{\xi}_I H_I$, then the image of $\xi$ under the operator $\partial_j$ is
$ \partial_j(\xi) = \sum_{I \in \mathcal{I}_N} \hat{\xi}_I H_{I-\delta_j},$
with the convention that if the multi-index $I-\delta_j$ has some negative
entries then $H_{I-\delta_j}=0$. For a function $f \in L^2(\RR^{\ZZ}, d\nu_0^{gc})$
 denote by $||f||=\sqrt{<f,f>}$ the $L^2$ norm of $f$.

The operators $\partial_j$ and $\partial_j^*$ act on the degree $N$ subspaces as follows:
 $$\partial_j ({\mathcal H}_N) \subseteq {\mathcal H}_{N-1}, \quad N \geq 1, \quad 
   \partial_j^* ({\mathcal H}_N) \subseteq {\mathcal H}_{N+1}, \quad N \geq 0.$$

The boundedness of these operators follows from the observation that 
$$\frac{1}{(N!)^N} \leq  \inf_{I \in \mathcal{I}_N} ||H_I||^2 \leq 
\sup_{I \in \mathcal{I}_N} ||H_I||^2 \leq 1, $$ 
and from the existence of two strictly positive constants, $C_1^N$, $C_2^N$, 
that depend just on $N$ such that
\begin{equation} \label{ineq1}
C_1^N \sum_{I \in \mathcal{I}_N} \hat{\xi}_I^2 \leq ||\xi||^2  \leq
C_2^N \sum_{I \in \mathcal{I}_N} \hat{\xi}_I^2, \quad \quad
C_1^N \sum_{I \in \mathcal{I}_N} \hat{\xi}_I^2 \leq ||\partial_j (\xi)||^2 \leq 
C_2^N \sum_{I \in \mathcal{I}_N} \hat{\xi}_I^2 .
\end{equation}
Indeed 
$$||\partial_j \xi ||^2 = \sum_{I \in \mathcal{I}_N} \hat{\xi}^2_I ||H_{I -\delta_j}||^2 \leq 
\sum_{I \in \mathcal{I}_N} \hat{\xi}_I^2 \leq (N!)^N \sum_{I \in \mathcal{I}_N}\hat{\xi}_I^2 ||H_I||^2 \leq (N!)^N ||\xi||^2.$$
and hence the norm of the operator $\partial_j :\mathcal{H}_N \to \mathcal{H}_{N-1}$
is bounded above by $(N!)^N$.

The vector field $D_0$ with constant coefficients has similar properties:
 $$D_0({\mathcal H}_N) \subseteq {\mathcal H}_{N-1}, \quad N \geq 1, \quad D_0^* 
  ({\mathcal H}_N) \subseteq {\mathcal H}_{N+1}, \quad N \geq 0.$$

For any function $\xi \in L^2(\RR^{\ZZ}, d\nu_0^{gc})$ and any test function 
$\phi \in {\mathcal H}_{N-1}$ we have
$$<D_n(\tau^m \xi ), \phi > = <\xi,\tau^{-m} (D_n^* \phi)> = 
 <\mathrm{Proj}_{{\mathcal H}_N}\xi, \tau^{-m} (D_n^* \phi)> = $$ 
\begin{equation} \label{qwe111}
= <D_n(\tau^m \mathrm{Proj}_{{\mathcal H}_N}\xi), \phi >, 
\end{equation}
$$<D_m(\tau^n \xi ), \phi > = <\xi, \tau^{-n} (D_m^* \phi)> = 
<\mathrm{Proj}_{{\mathcal H}_N}\xi, \tau^{-n} (D_m^* \phi)> = $$ 
\begin{equation} \label{qwe1}
= <D_m(\tau^n \mathrm{Proj}_{{\mathcal H}_N}\xi), \phi >.
\end{equation}
It follows that $D_n(\tau^m \xi ) =D_m(\tau^n \xi )$ in the weak sense if and only if 
$D_n(\tau^m \mathrm{Proj}_{{\mathcal H}_N}\xi)=D_m(\tau^n \mathrm{Proj}_{{\mathcal H}_N}\xi)$ 
in the strong sense for all $N \geq 0$.

We recall that a function $\xi$ is closed if and only if  
$D_n(\tau^m \xi ) =D_m(\tau^n \xi )$ for all $m,n \in \ZZ$, which, by the 
previous equalities (\ref{qwe111}) and (\ref{qwe1}), is equivalent to 
$$D_n(\tau^m \mathrm{Proj}_{{\mathcal H}_N}\xi)=D_m(\tau^n \mathrm{Proj}_{{\mathcal H}_N}\xi)
 \quad m,n \in \ZZ \quad  N \geq 0.$$ 
Therefore, a function $\xi$ is closed if 
and only if $\mathrm{Proj}_{{\mathcal H}_N}\xi$ is closed for all $N \geq 0$. \hfill \qed

{\it Note.} If $\xi =\sum_{I \in \mathcal{I}_N} \hat{\xi}_I H_I$ is a function inside the space 
$\mathcal{H}_N$, two norms can be defined for $\xi$: the $L^2$ norm $||\xi||$ and 
the sum of squared Fourier coefficients $\sum_{I \in \mathcal{I}_N} \hat{\xi}_N^2$.
It is important to note  the inequality (\ref{ineq1}) implies that these two norms define 
the same topology on the space $\mathcal{H}_N$.

{\it Note.} Assume $\xi \in \mathcal{H}_N$ is a $D_0$-closed function, with Fourier series expansion 
$\xi = \sum_{I \in \mathcal{I}} \hat{\xi}_I H_I$. We calculate,
$$
D_n \xi  =   \sum_{I \in \cI} \bigg[ \sum_{k \in \ZZ} a_k \hat{\xi}_{I+\delta_{(n+k)}}\bigg] H_I,  \quad
D_0(\tau^n \xi)  =  \sum_{I \in \cI} \bigg[ \sum_{k \in \ZZ} a_k \hat{\xi}_{\tau^n(I+\delta_{k})} \bigg] H_I. 
$$
Therefore a function is closed if and only if its Fourier coefficients satisfy 
the relations:
\begin{equation} \label{feq4}
\sum_{k \in \ZZ} a_k \hat{\xi}_{I+\delta_{(n+k)}}=\sum_{k \in \ZZ} a_k \hat{\xi}_{\tau^n(I+\delta_{k})}
\quad n \in \ZZ, \quad I \in {\mathcal I}.
\end{equation}

{\bf Construction of exact functions.} It is important to have some examples 
of functions that are exact. The functions that will be constructed next 
will be used in the proof of the decomposition theorem \ref{the1}, to 
approximate closed functions with exact ones. 
\begin{lemma} \label{lem2}
Let $c$ be a function defined on the set of orbits with finite support 
(i.e., c(o)=0 except for finitely many orbits $o$). The function 
$$ \xi = \sum_{o \in \mathcal{O}} c(o ) D_0 \bigg[ \sum_{n \in \ZZ} \tau^n 
H_{R(o)+\delta_0} \bigg]$$ is $D_0$-exact and the Fourier coefficients of $\xi$ are
\begin{equation} \label{feq85}
\hat{\xi}_I= \sum_{k \in \ZZ} a_k c(\co{\tau^{-k} I}).
\end{equation}
\end{lemma}

{\bf Proof.} The function $\xi$, that has been introduced is well-defined
since the sum is over a finite set, and is exact as a sum of 
exact functions. To conclude the lemma we need to calculate the Fourier 
coefficients of $\xi$. We have, 
\begin{eqnarray} \label{feq2}
\xi & = & \sum_{o \in \mathcal{O}} c(o) D_0\bigg[ \sum_{n \in \ZZ} 
         H_{ \tau^{-n}(R(o)+\delta_0)} \bigg] 
     =  \sum_{o \in \mathcal{O} , n \in \ZZ} c(o) \sum_{k \in \ZZ}
          a_kH_{\tau^{-n}(R(o)+\delta_0- \delta_{-n+k})} = \nonumber \\
    &  & = \sum_{o \in \mathcal{O} , n \in \ZZ} c(o) 
          \sum_{k \in \ZZ} a_kH_{\tau^{-k}[(-n+k)\cdot R(o)]}  
    = \sum_{I \in \mathcal{I}} \sum_{k \in \ZZ} a_kc(\co{\tau^{k}I})H_I.
%   & = & \sum_{o \in \mathcal{O}}c(o) \bigg[ 
%          \sum_{n-1 \in s(R(o))} H_{\tau[(n-1) \cdot R(o)]}
%         -2\sum_{n-1 \in s(R(o))} H_{n\cdot R(o)} 
%          + \sum_{n+1 \in s(R(o))} H_{\tau^{-1}[(n+1) \cdot R(o)]}  \bigg].
\end{eqnarray} 

To justify the integration by parts in the calculation above (\ref{feq2}) 
we make the following observation. For any multi-index $I \in \mathcal{I}$ 
there exists a unique orbit $o \in \mathcal{O}$ and a unique integer $n\in \ZZ$ 
such that $I=\tau^{-k}[(-n+k) \cdot R(o)]$. This is a consequence of the freeness of 
the action (\ref{action}). Moreover, the orbit $o$ is the same as $\co{\tau^k I}$. 
We stress again that the sums in (\ref{feq2}) are over finite sets as $c$ has 
finite support. Actually all computations that we carried out 
to prove this lemma are valid because $c$ is a function with finite support 
and the sums are finite, although this wasn't emphasized each time we used it.
Also we have made use of the convention that $H_I=0$ if $I$ is a multi-index
with negative entries. \hfill \qed

\begin{lemma}
Let $N \geq 1$ be a natural number and $e=(1, \dots , 1) \in \ZZ^N$.  In addition if $\ww{c}$ is 
a real-valued function defined on $\ZZ^N$, with finite support and $\ww{S}_N$-invariant then 
the function 
\begin{equation} \label{exact}
\xi_{\ww{c}} = \sum_{I \in \mathcal{I}_N} \bigg(\sum_{k \in \ZZ}a_k \ww{c}(z_I-ke)\bigg) H_I
\end{equation}
is a well-defined $D_0$-exact function in the degree $N$ subspace $\mathcal{H}_N$.
\end{lemma}

\textbf{Proof.} This lemma follows from lemma \ref{lem2}.  Since $\ww{c}:\ZZ^N \to \RR$ 
is $\ww{S}_N$-invariant, it makes sense to introduce $c:\mathcal{O} \to \RR$, where 
$c(o)= \ww{c}(z_I)$ if $I$ is a multi-index in the orbit $o$ of degree $N$, and $c(o)=0$ 
otherwise. We should note that if $I$ is a multi-index
in the orbit $o$ then $z_I+ke = z_{\tau^{-k}I}$ and $\ww{c}(z_I-ke)=c( o(\tau^{k}I))$.  
Hence the Fourier coefficients of 
the function $\xi_{\ww{c}}$ are of the form $\sum_{k \in \ZZ}a_k c(o(\tau^{k}I))$, and 
the function $\xi_{\ww{c}}$ is $D_0$-exact. \hfill \qed

In the previous lemma an operator has come out in a natural way in our 
construction of exact functions. Below we provide the exact definition 
of this operator. 
\begin{definition}
Let $D_0=\sum_{k \in \ZZ} a_k \partial_k$ be a vector field with constant 
coefficients, all the coefficients being zero except finitely many.
The vector field $D_0$ defines an operator $T_{D_0}$ that acts 
on functions $c:\ZZ^N \to \RR$ and produces a function 
$T_{D_0}c:\ZZ^N \to \RR$, where
$$(T_{D_0}c)(z) = \sum_{k \in \ZZ} a_k c(z-ke), \quad z \in \ZZ^N.$$
Above $e$ is the vector $(1,\dots ,1)$ of the lattice $\ZZ^N$.
\end{definition}

%%%%%%%%%%%%%%%%%%%%%%%%%%%%%%%%%
%\section{Proof of the decomposition theorem}
%%%%%%%%%%%%%%%%%%%%%%%%%%%%%%%%%
\section{Proof of the decomposition theorem \ref{the1}}

We start by listing two important properties of the operator $T_{D_0}$
introduced at the end of the previous section.

\begin{lemma} \label{compact}
Let $c$ be a real-valued function defined on the lattice $\ZZ^N$, $N\geq 1$.
We assume that the function $c$ is square-summable and $\widetilde{S}_N$- 
invariant.  Then, there exists a sequence $(c_n)_{n \geq 1}$ of real-valued, 
finitely supported, $\widetilde{S}_N$-invariant functions such that
$T_{D_0}c_n \to T_{D_0}c$ as $n \to \infty$ and the convergence is in the Hilbert space
topology of $L^2(\ZZ^N)$. 
\end{lemma}

{\bf Proof.} We define a sequence of $\widetilde{S}_N$-invariant 
regions of the lattice $\ZZ^N$, namely
\begin{equation} \label{hexa}
P_i = \cup_{\gamma \in \widetilde{S}_N} \gamma\{z=(z_1, \dots , z_N) \in \ZZ^N | 0 \leq z_1 \leq \dots \leq z_N \leq i-1 \},  \quad i\geq 1.
\end{equation} 
For the reader convenience we add  two pictures of the region $P_i$ in dimension $N=1$, respectively $N=2$.
In dimension $N=1$ the region $P_i$ contains the lattice points inside the segment $[-i+1,i-1]$, whereas
in dimension $N=2$ the region $P_i$ contains the lattice points inside the hexagon shown below.

\begin{center}
\begin{picture}(400,50)(-50,0)
%\put(10,30){\circle*{2}}
%\put(30,30){\circle*{2}}
%\put(50,30){\circle*{2}}
\put(70,30){\circle*{5}}
\put(90,30){\circle*{5}}
\put(110,30){\circle*{5}}
\put(130,30){\circle*{5}}
\put(150,30){\circle*{5}}
\put(170,30){\circle*{5}}
\put(190,30){\circle*{5}}
\put(210,30){\circle*{5}}
\put(230,30){\circle*{5}}
%\put(250,30){\circle*{2}}
%\put(270,30){\circle*{2}}
%\put(290,30){\circle*{2}}
\put(60,10){-i+1}
\put(82,10){-i+2}
\put(120,10){$\dots$}
\put(150,10){0}
\put(175,10){$\dots$}
\put(205,10){i-2}
\put(225,10){i-1}
\end{picture}
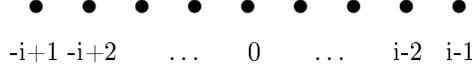
Figure 1. The region $P_i$ in dimension $N=1$.
\end{center}

\begin{center}
\begin{picture}(400,250)(-50,0)
\put(70,30){\circle*{5}}
\put(90,30){\circle*{5}}
\put(110,30){\circle*{5}}
\put(130,30){\circle*{5}}
\put(150,30){\circle*{5}}
\put(70,50){\circle*{5}}
\put(70,70){\circle*{5}}
\put(70,90){\circle*{5}}
\put(70,110){\circle*{5}}
\put(150,190){\circle*{5}}
\put(170,190){\circle*{5}}
\put(190,190){\circle*{5}}
\put(210,190){\circle*{5}}
\put(230,190){\circle*{5}}
\put(230,170){\circle*{5}}
\put(230,150){\circle*{5}}
\put(230,130){\circle*{5}}
\put(230,110){\circle*{5}}
\put(90,130){\circle*{5}}
\put(110,150){\circle*{5}}
\put(130,170){\circle*{5}}
\put(170,50){\circle*{5}}
\put(190,70){\circle*{5}}
\put(210,90){\circle*{5}}
\put(150,170){\circle*{5}}\put(170,170){\circle*{5}}\put(190,170){\circle*{5}}\put(210,170){\circle*{5}}
\put(130,150){\circle*{5}}\put(150,150){\circle*{5}}\put(170,150){\circle*{5}}\put(190,150){\circle*{5}}
\put(210,150){\circle*{5}}
\put(110,130){\circle*{5}}\put(130,130){\circle*{5}} \put(150,130){\circle*{5}}\put(170,130){\circle*{5}} 
\put(190,130){\circle*{5}}\put(210,130){\circle*{5}}
\put(90,110){\circle*{5}}\put(110,110){\circle*{5}} \put(130,110){\circle*{5}}\put(150,110){\circle*{5}} 
\put(170,110){\circle*{5}}\put(190,110){\circle*{5}} \put(210,110){\circle*{5}}
\put(90,90){\circle*{5}}\put(110,90){\circle*{5}} \put(130,90){\circle*{5}}\put(150,90){\circle*{5}} 
\put(170,90){\circle*{5}}\put(190,90){\circle*{5}} 
\put(90,70){\circle*{5}}\put(110,70){\circle*{5}} \put(130,70){\circle*{5}}\put(150,70){\circle*{5}} 
\put(170,70){\circle*{5}}\put(90,50){\circle*{5}}\put(110,50){\circle*{5}}\put(130,50){\circle*{5}} 
\put(150,50){\circle*{5}} 
\put(50,15){(-i+1,-i+1)}
\put(50,115){(-i+1,0)}
\put(140,200){(0,-i+1)}
\put(220,200){(i-1,i-1)}
\put(220,115){(i-1,0)}
\put(140,15){(0,-i+1)}
\put(140,115){(0,0)}
%\put(50,10){\circle*{2}} \put(50,30){\circle*{2}}\put(50,50){\circle*{2}}
%\put(50,70){\circle*{2}} \put(50,90){\circle*{2}}\put(50,110){\circle*{2}}
%\put(70,130){\circle*{2}} \put(90,150){\circle*{2}}\put(110,170){\circle*{2}}
%\put(130,190){\circle*{2}} \put(150,210){\circle*{2}}\put(170,210){\circle*{2}}
%\put(190,210){\circle*{2}} \put(210,210){\circle*{2}}\put(230,210){\circle*{2}}\put(250,190){\circle*{2}}
%\put(250,170){\circle*{2}} \put(250,150){\circle*{2}}\put(250,130){\circle*{2}}
%\put(250,210){\circle*{2}}
\end{picture}
Figure 2. The region $P_i$ in dimension $N=2$.
\end{center}

Beside being $\ww{S}_N$-invariant, the sequence of regions $(P_i)_{i \geq 1}$ defined above, grows 
to cover the entire lattice $\ZZ^N$ as $i \to \infty$. 
Define $c_n$ to be $c {\bf 1}_{P_n}$, for $n \geq 1$.
Since ${\bf 1}_{P_n}$ is the characteristic function of the region $P_n$ we have 
immediately that $c_n$ is a finitely supported, $\ww{S}_N$-invariant function. 
Square-summability of $c$ implies that $c_n \to c$ as $n \to \infty$ in the topology 
of $L^2(\ZZ^N)$ (the norm $||c-c_n||^2=\sum_{z \notin P_n} c^2(z)$ involves only the values of $c$ outside 
the region $P_n$, and these values decay to zero as $n \to \infty$ since $c$ is 
square-summable). Then, obviously, $c_n \to c$ and $Tc_n \to Tc$ as $n \to \infty$  
in the topology of $L^2(\ZZ^N)$.  \hfill \qed

Below we discuss certain facts about the Fourier transform of functions 
defined on the lattice $\ZZ^N$.
The Fourier transform of a function $c:\ZZ^N \to \mathbb{C}$ is formally defined to be
$$ \mathcal{F} c :[-\pi , \pi )^N \to \mathbb{C}, \quad \mathcal{F}c (\alpha ) = 
\frac{1}{\sqrt{2\pi }}\sum_{z \in \ZZ^N} c(z) e^{i z \alpha }.$$
In the exponent above $z\alpha$ stands for the dot product $z_1\alpha_1+\dots +z_N \alpha_N$.
The reader may consult Rudin \cite{Rud} for an extended treatment of Fourier transform
of functions defined on lattice. We remind the reader that
$\mathcal{F}$ is an isometry between the spaces  $L^2( \ZZ^N )$ and $L^2 ([-\pi, \pi )^N)$.
%$ || c ||_{L^2(\ZZ^N)} = || \mathcal{F} (c) ||_{L^2 [ -\pi , \pi )^N}$
 %$\mathcal{F} (T_j c) (\alpha ) = e^{-i(\alpha_1+\dots +\alpha_N)} \mathcal{F} (c)(\alpha)$.
The space $L^2([-\pi,\pi)^N)$ is considered with respect to the Lebesgue 
measure on $[-\pi, \pi)^N$.
Also if $c$ is invariant under a certain group of transformations then $\mathcal{F}c$ 
is invariant, as well. Though the symmetry group of $\mathcal{F}c$ might not 
coincide with the symmetry group of $c$. Indeed if $c$ is symmetric, or $S_N$-invariant
then $\mathcal{F}c$ is symmetric. Now suppose that $c$ is $\ww{S}_N$-invariant
then $\mathcal{F}c$ is invariant under the action of the group $\ww{\Sigma}_N$,
 generated by the transformations:
$$s_{ii+1}:[-\pi,\pi)^N \to [-\pi,\pi)^N, \quad 1 \leq i \leq N-1$$
$$s_{ii+1}(\alpha_1, \dots , \alpha_N) =(\alpha_1, \dots , \alpha_{i+1}, \alpha_{i}, \dots , \alpha_N), $$
and 
$$g:[-\pi,\pi)^N \to [-\pi,\pi)^N,
g(\alpha_1, \dots , \alpha_N) =(\mathrm{mod}_{2\pi}(-\alpha_1- \dots -\alpha_N) , \alpha_{2}, \dots , \alpha_N).$$
On the line above we used the notation $\mathrm{mod}_{2\pi}(t)$. Any real number $t$ can be written 
uniquely as $2\pi a+b$, where $a$ is an integer number and $b$ is a real number in 
the interval $[-\pi, \pi)$. By $\mathrm{mod}_{2\pi}(t)$ we denote the remainder $b$.
It is also true that if the Fourier transform $\mathcal{F}c$ is $\ww{\Sigma}_N$-invariant
then $c$ is $\ww{S}_N$-invariant. 
 
In section \ref{clos-exac} we have established that a function 
$\xi=\sum_{I \in \mathcal{I}_N} \hat{\xi}_I H_I$ is $D_0$-closed if and only if the following holds:
\begin{equation} \label{feq10}
\sum_{k \in \ZZ} a_k \hat{\xi}_{I+\delta_{(n+k)}}=\sum_{k \in \ZZ} a_k \hat{\xi}_{\tau^n(I+\delta_{k})}
\quad n \in \ZZ, \quad I \in {\mathcal I}.
\end{equation}
Obviously we can use the Fourier coefficients of $\xi$ to construct a $S_N$-invariant
function $\hat{\xi}: \ZZ^N \to \RR$, $\hat{\xi}(z)=\hat{\xi}_I$ if $z \stackrel{S_N}{\sim} z_I$. 
The relations (\ref{feq10}) force our function $\xi$ to satisfy
\begin{equation} \label{feq11}
\sum_{k \in \ZZ} a_k \hat{\xi}(z+ke_1) = \sum_{k \in \ZZ} a_k \hat{\xi}(z-z_1e - (z_1+k)e_1), 
\quad z=(z_1, \dots , z_N) \in \ZZ^N.
\end{equation}
The vectors $e$ and $e_1$ of the lattice $\ZZ^N$ are $(1, \dots , 1)$ and $(1, 0, \dots , 0)$, respectively.
After applying the Fourier transform in both sides of the equation (\ref{feq11}) 
we find that $\hat{\xi}$  satisfies 
\begin{equation} \label{feq12}
p(e^{-i\alpha_1}) (\mathcal{F}\hat{\xi})(\alpha) = p(e^{i(\alpha_1+\dots + \alpha_N)})(\mathcal{F}\hat{\xi})(g( \alpha)) , \quad 
\alpha=(\alpha_1, \dots , \alpha_N) \in [-\pi, \pi)^N,
\end{equation}
where $p$ is the rational function $p(x)=\sum_{k \in \ZZ} a_k x^k$.
To wrap up our argument we can say that the $D_0$-closedness condition implies property (\ref{feq12}).

Next we shall establish a crucial fact about  functions that satisfy relation (\ref{feq12}).
\begin{lemma}
Let $\hat{\xi}$ be a real-valued, $S_N$-invariant function defined on the lattice 
$\ZZ^N$, $N \geq 1$ that satisfies (\ref{feq12}). Then, 
there exists a sequence $(c_n)_{n \geq 1}$ of real-valued, square-summable, 
$\ww{S}_N$-invariant functions defined on the lattice $\ZZ^N$ such that 
$T_{D_0}c_n \to \hat{\xi}$ as $n \to \infty$ in the topology of $L^2(\ZZ^N)$.
\end{lemma}

{\bf Proof.} 
Examples of functions $\hat{\xi}$ satisfying the properties listed in the hypothesis
of this lemma, are functions constructed, as explained before in this section,
 from the Fourier coefficients of closed functions.

The first step towards establishing our result is to solve, at least on a formal level,
 the equation $T_{D_0}c =\hat{\xi}$. 
The Fourier transform for functions 
defined on the lattice will help us to make our guess.

After applying the Fourier transform in each side of the equation $T_{D_0}c =\hat{\xi}$, we get 
$$p(e^{i(\alpha_1+\dots +\alpha_N)})(\mathcal{F}c)(\alpha) =(\mathcal{F}\hat{\xi})(\alpha), \quad
\alpha=(\alpha_1, \dots , \alpha_N) \in \ZZ^N .$$
where $p$ is the rational function $\sum_{j \in \ZZ} a_j x^j$ canonically associated to $T_{D_0}$.

After multiplying with a high enough power of $x$ the equation $p(x)=0$, we see that any 
solution $x$ of $p(x)=0$ has to be the root of a certain polynomial. Since 
the number of roots of any polynomial is finite, the number of solutions of $p(x)=0$ 
is finite, as well.  
If the equation $p(x)=0$ has no solutions on the unit circle then the equation $T_{D_0}c=\hat{\xi}$ 
can be solved in the space $L^2(\ZZ^N)$. Indeed the unique, square-summable
solution of $T_{D_0}c=\hat{\xi}$ is 
$$c=\mathcal{F}^{-1}\bigg( \frac{1}{p(e^{i(\alpha_1+\dots +\alpha_N)})}( \mathcal{F}\hat{\xi}) \bigg).$$
That $\hat{\xi}$ has been built out of the Fourier coefficients of a closed 
function implies that $c$ is $\ww{S}_N$-invariant. Hence the lemma is proved in this 
case. We can choose $c_n=c$, for any $n \geq 1$.

A more involved case is when the equation $p(x)=0$ has solutions on the unit circle. 
If the sum of the coefficients of $p$ is equal to $0$ then the number $1$ is a 
solution of  $p(x)=0$. The cases arising from interacting particle models are of this kind. 
The difficulty in this case arises because the 
equation $Tc=\hat{\xi}$ can not be solved in $L^2(\ZZ^N)$. To get around this 
problem we shall consider a slightly modified equation $\mathcal{F}(T_{D_0}c)=(\mathcal{F}\hat{\xi})1_{A_n}$.
The function $\mathcal{F}\hat{\xi}$ is multiplied with the characteristic function of a 
set $A_n$ that in $\ww{\Sigma}_N$-invariant and carefully chosen to avoid the unit roots of $p$. More precisely,
$$A_n= \cap_{\gamma \in \ww{\Sigma}_N} \bigg\{ \gamma(\alpha)\; \bigg|\; \alpha =(\alpha_1, \dots , \alpha_N) \in [-\pi,\pi]^N, $$ 
$$| \; \mathrm{mod}_{2\pi}(\alpha_1+\dots +\alpha_N)-r_k|> \frac{1}{n}, \;\; e^{ir_k} 
      \mathrm{unit}\; \mathrm{root} \; \mathrm{of}\; p \bigg\}.$$  

The next table contains the roots of the rational function $p(x)$ in four particular cases. \\
\begin{tabular}{l|l|l}
Vector field $D_0$ & Rational function $p(x)$ & Solutions of $p(x)=0$ \\ \hline
$\partial_0$            & $1$    & none  \\
$\partial_1 - \partial_0$ & $x-1$ & $1$ \\
$\partial_1-2\partial_0+\partial_{-1}$ & $x-2+x^{-1}$ & $1$, $1$\\
$\partial_3-\partial_0$ & $x^3-1$ & $1$, $\frac{1+\sqrt{-3}}{2}$, $\frac{1-\sqrt{-3}}{2}$ 
\end{tabular}

If the vector field $D_0$ is $\partial_1 -\partial_0$ or $\partial_1-2 \partial_0+\partial_{-1}$ 
then the region $A_n$ in dimension $N=1$ is just
$A_n= \{ \alpha \in [-\pi,\pi) \; | \; |\alpha| > \frac{1}{n} \}$, \\
%$ N=2, \quad A_n= \{ \alpha =(\alpha_1, \alpha_2) \in [-\pi,\pi]^2 \; | \; |\alpha_1+\alpha_2| > \frac{1}{n} \}$. We include a picture of the region $A_n$ in dimension $N=1$ and respectively $N=2$.

Suppose we are given a $Y_0$-closed function $\xi$ with Fourier expansion 
$\xi= \sum_{I \in \mathcal{I}}\hat{\xi}_I H_I$. We can actually turn our 
graph into a weighted graph by assigning to each directed edge $(o(I),o(\tau I))$
the weight $\hat{\xi}_I$.  

\begin{center}
\begin{picture}(400,50)(-50,0)
\put(50,27){[}
\put(130,27){)}
\put(170,27){(}
\put(250,27){]}
\put(50,30){\line(1,0){80}}
\put(170,30){\line(1,0){80}}
\put(50,10){-$\pi$}
\put(150,10){0}
\put(120,10){$-\frac{1}{n}$}
\put(170,10){$\frac{1}{n}$}
\put(250,10){$\pi$}
\end{picture}
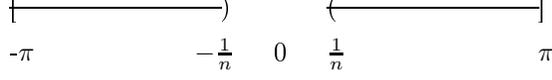
Figure 3. The region $A_n$ in dimension $N=1$.
\end{center}

%\begin{center}
%\begin{picture}(400,250)(-50,0)
%\put(50,30){\line(1,0){200}}
%\put(50,30){\line(0,1){200}}
%\put(50,230){\line(1,0){200}}
%\put(250,30){\line(0,1){200}}
%\put(230,30){\line(-1,1){180}}
%\put(250,50){\line(-1,1){180}}
%\put(110,30){\circle*{5}}
%\put(130,30){\circle*{5}}
%\put(250,210){\line(-1,1){20}}
%\put(70,30){\line(-1,1){20}}
%\end{picture}
%Figure 4. The region $A_n$ in dimension $N=2$.
%\end{center}

\begin{center}
\begin{picture}(400,250)(-50,0)
%\linethickness{4pt}
\put(70,30){\line(1,0){60}}
\put(130,30){\line(0,1){80}}
\put(50,50){\line(0,1){60}}
\put(50,110){\line(1,0){80}}
\put(70,30){\line(-1,1){20}}

\put(170,30){\line(1,0){60}}
\put(170,30){\line(0,1){60}}
\put(230,30){\line(-1,1){60}}

\put(250,50){\line(-1,1){60}}
\put(250,50){\line(0,1){60}}
\put(250,110){\line(-1,0){60}}

\put(250,210){\line(0,-1){60}}
\put(230,230){\line(-1,0){60}}
\put(250,210){\line(-1,1){20}}
\put(250,150){\line(-1,0){80}}
\put(170,150){\line(0,1){80}}

\put(50,150){\line(1,0){60}}
\put(50,150){\line(0,1){60}}
\put(110,150){\line(-1,1){60}}

\put(130,170){\line(-1,1){60}}
\put(130,230){\line(0,-1){60}}
\put(130,230){\line(-1,0){60}}

\end{picture}
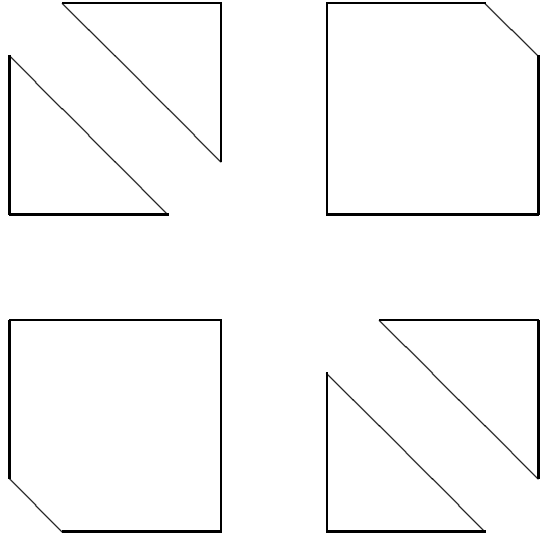
Figure 4. The region $A_n$ in dimension $N=2$, as a subset of the square $[-\pi,\pi]^2$.
\end{center}

%\begin{center}
%\begin{figure} 
%\begin{psfigure}(400,250)(-50,0)
%\pspolygon*(250,50)(250,110)
%\end{psfigure}
%\caption{The region}
%\end{figure}
%\end{center}

Let $c_n$ be the unique $L^2(\ZZ^N)$ solution of the equation 
$\mathcal{F}(T_{D_0}c_n)=(\mathcal{F}\hat{\xi})1_{A_n}$, $n \geq 1$. 
The solution $c_n$ is defined through
$$ c_n = \mathcal{F}^{-1}\bigg( \frac{1}{p(e^{i(\alpha_1+\dots +\alpha_N)})} 
               (\mathcal{F}\hat{\xi})(\alpha)1_{A_n}(\alpha) \bigg).$$
The $\ww{\Sigma}_N$-invariance of the set $A_n$ and the fact that $\hat{\xi}$ is constructed 
from the Fourier coefficients of a closed function implies that $c_n$ is $\ww{S}_N$-invariant, 
$n\geq 1$.

Obviously we have the convergence $\mathcal{F}(T_{D_0}c_n) \to \mathcal{F}\hat{\xi}$ as 
$n \to \infty$ in the topology of $L^2([-\pi,\pi]^N)$, hence $T_{D_0}c_n \to \hat{\xi}$ as 
$n \to \infty$ in the topology of $L^2(\ZZ^N)$.  \hfill \qed

\textbf{Proof of the decomposition theorem \ref{the1}}. 
Let $\xi \in L^2(\RR^{\ZZ}, d\nu_0^{gc})$ be a $D_0$-closed function. From lemma
4.1 we know that $\mathrm{Proj}_{\mathcal{H}_N}\xi$  is $D_0$-closed, for any $N \geq 0$.
$\mathrm{Proj}_{\mathcal{H}_0}\xi$ is a constant function. If the sum of the coefficients 
of $D_0$ is not equal to zero then $\mathrm{Proj}_{\mathcal{H}_0}\xi$ is $D_0$-exact, 
otherwise $\mathrm{Proj}_{\mathcal{H}_0}\xi$ is orthogonal on any $D_0$-exact function. 
Therefore, the  decomposition theorem follows as long as we establish that any $D_0$-closed
function $\xi \in \mathcal{H}_N$ can be approximated by $D_0$-exact functions, for any $N \geq 1$.

Assume $\xi =\sum_{I \in \mathcal{I}_N} \hat{\xi}_I H_I \in \mathcal{H}_N$, $N \geq 1$. 
Define the $S_N$-invariant function $\hat{\xi} :\ZZ^N \to \RR$ through $\hat{\xi}(z) = \hat{\xi}_I$ if 
$z \stackrel{S_N}{\sim} z_I$, see (3). We use lemmas 5.1 and 5.2 to find a 
sequence of finitely supported, $\ww{S}_N$-invariant functions $(c_n)_{n \geq 1}$, 
such that $T_{D_0} c_n  \to \hat{\xi}$ as $n \to \infty$ in the topology of $L^2(\ZZ^N)$.
But lemmas 4.2 and 4.3 tell us that each of $T_{D_0}c_n$, $n \geq 1$ defines a 
$D_0$-exact function $\xi_{c_n}$, see (\ref{exact}). At the end of lemma 4.1 we noticed
that the topology of $\mathcal{H}_N$ and $L^2(\ZZ^N)$ are equivalent, hence 
we can claim that $\xi_{c_n} \to \xi$ as $n \to \infty$ in the Hilbert space topology
of $\mathcal{H}_N$, or $L^2(\RR^{\ZZ}, d\nu_0^{gc})$. \hfill \qed

%%%%%%%%%%%%%%%%%%%%%%%%%%%%%%%%%%%%%%%
%\section{Second-order Ginzburg-Landau field and algebraic topology}
%%%%%%%%%%%%%%%%%%%%%%%%%%%%%%%%%%%%%%%
\section{Second-order Ginzburg-Landau field and algebraic topology}

We conclude with some remarks about the second-order 
Ginzburg-Landau field $Y_0=\partial_1-\partial_0$ which has been studied in the work 
of S. R. S. Varadhan, \cite{Var1}. Our approach places Varadhan's result 
in a new light  by depicting an algebraic topologic aspect, to be 
explained below.

In section \ref{multiind} we presented an extensive study of the set of multi-indices $\mathcal{I}$. 
There we partitioned the set of multi-indices $\mathcal{I}$ into disjoint orbits, and 
we denoted by $\mathcal{O}$ the space of orbits. Below we exhibit a procedure 
to construct  a directed graph that has as vertices the orbits of the set of multi-indices.

A directed graph is a pair $(V, E)$ of two sets, where $V$ is the set of vertices 
of the graph and $E$ is the set of directed edges. A directed edge is 
a pair of two vertices $(v_1, v_2)$ where the first vertex indicates the starting point 
of the edge and  the second vertex indicates the tip of the edge. For us we choose $V$ 
to be the set of orbits $\mathcal{O}$. Also we say that we have a directed 
edge $(o_1, o_2)$ if there exist  a multi-index $I \in o_1$ such that $\tau I \in o_2$. 
Notice that if there exists an edge between two orbits then the orbits contain 
multi-indices with
identical degrees. Hence our graph will have at least one connected component
for each degree $N \geq 0$. We shall show that there exists precisely one 
connected component for each degree $N \geq 0$. 
 
We would like 
to have a concrete or geometric presentation of the graph. For this purpose 
we use the identification of the set of orbits $\mathcal{O}_N$ containing the 
multi-indices of  degree $N$, with the cone $\mathcal{C}_N^+$ of the lattice 
$\ZZ^N$, $N \geq 0$.  

Assume $N=0$, then 
$\mathcal{O}_N$ contains a single orbit, the orbit of the 
multi-index $0$ and this orbit contains a single multi-index. 
Since the multi-index $0$ has the property that $0=\tau 0$, 
we will have a directed edge going out of and returning to $0$; 
in other words we have a loop at $0$. 

Assume that $N \geq 1$. 
It can be shown that a directed edge links $z \in \mathcal{C}_N^+$ to 
$z' \in \mathcal{C}_N^+$ if and only if either $z'=z-e_1-\dots -e_N$ 
or $z'=z+e_i$ for some $1 \leq i \leq N$. Here $e_i$ is lattice 
vector $(0,\dots , 1, \dots, 0)$ with the $i$th coordinate $1$.   
We shall indicate below how this presentation of the graph can be obtained.

Let $o \in \mathcal{O}_N$ be some orbit and $R(o) =\sum_{i=1}^k a_i \delta_{n_i}$, 
with $a_i \geq 1$ and $0 \leq n_1 < n_2< \dots < n_k$ be the 
representative of the orbit $o$.  Given our rule, the orbit $o$  is connected
by an edge going out of $o$ to each of the orbits $ o(\tau R(o))$,  $o(\tau(n_1 \cdot R(o)))$,  
$\dots$,  $o(\tau(n_k \cdot R(o)))$.  For each of the orbits in the list before
we can calculate the representatives and the corresponding 
point in the cone $\mathcal{C}_N^+$.

For example: the representative of the orbit $o$ is 
$R(o) =\sum_{i=1}^k a_i \delta_{n_i}$ and the cone point is 
$$z_{R(o)}=(\underbrace{n_1, \dots , n_1}_{a_1}, \dots , \underbrace{n_k, \dots , n_k}_{a_k}).$$
Assume $n_1 \geq 1$. The representative of the orbit of $\tau R(o)$ 
is $ \tau R(o) $ and the corresponding cone point is 
 $$z_{\tau R(o)}=(\underbrace{n_1-1, \dots , n_1-1}_{a_1}, \dots , \underbrace{n_k-1, \dots , n_k-1}_{a_k}).$$
We notice that $z_{\tau R(o)} = z_{R(o)}-e_1-\dots - e_N$. 
Also if $n_1=0$ the representative of the orbit of $\tau R(o)$ 
is $ (-1) \cdot \tau R(o) $ and the corresponding cone point is 
 $$z_{(-1) \cdot \tau R(o)}=(\underbrace{0, \dots , 0}_{a_1-1}, 1, \dots , \underbrace{n_k-1, \dots , n_k-1}_{a_k}),$$
and  $z_{(-1) \cdot \tau R(o)} = z_{R(o)}+e_{a_1}$.  Similarly, we can analyze the other
orbits connected with $o$. 

In particular our discussion proves that for any two given orbits 
$o_1$ and $o_2$ if there exists a multi-index $I \in o_1$ and $\tau I \in o_2$ then this 
multi-index is unique. As we will see later that this observation allows us to assign 
in a unique way a multi-index to any directed edge of our graph.  

We include three pictures of the connected components of the  directed graph
for $N=0$, $N=1$ and $N=2$.

\begin{center}
\begin{picture}(400,50)(-50,0)
\put(150,60){\vector(1,0){5}}
\put(150,40){\circle{40}}
\put(150,20){\circle*{3}}
%\put(145,22){\vector(2,-1){7}}
\put(149,8){0}
\end{picture}
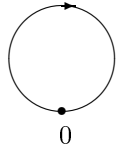
Figure 5. The connected component of the graph for $N=0$.
\end{center}

\begin{center}
\begin{picture}(400,50)(-50,0)
\put(170,30){\vector(1,0){40}}
\put(50,30){\vector(1,0){40}}
\put(90,30){\vector(1,0){40}}
\put(130,30){\vector(1,0){40}}
\qbezier(50,30)(70,80)(90,30)
\qbezier(90,30)(110,80)(130,30)
\qbezier(130,30)(150,80)(170,30)
\qbezier(170,30)(190,80)(210,30)
\put(50,30){\circle*{3}}
\put(90,30){\circle*{3}}
\put(130,30){\circle*{3}}
\put(170,30){\circle*{3}}
\put(210,30){\circle*{3}}

\put(220,30){\dots}
\put(50,10){0}
\put(90,10){1}
\put(130,10){2}
\put(170,10){3}
\put(210,10){4}

\put(70,55){\vector(-1,0){5}}
\put(110,55){\vector(-1,0){5}}
\put(150,55){\vector(-1,0){5}}
\put(190,55){\vector(-1,0){5}}
\end{picture}
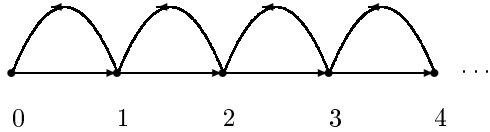
Figure 6. The connected component of the graph for $N=1$.
\end{center}

\begin{center}
\begin{picture}(400,200)(-50,0)
\put(170,30){\vector(1,0){40}} \put(210,30){\vector(0,1){40}} \put(210,70){\vector(-1,-1){40}} \put(210,70){\vector(0,1){40}}  
  \put(210,110){\vector(-1,-1){40}}  \put(210,110){\vector(0,1){40}}  \put(210,150){\vector(-1,-1){40}} \put(210,150){\vector(0,1){40}} 
  \put(210,190){\vector(-1,-1){40}}
\put(50,30){\vector(1,0){40}} \put(90,30){\vector(0,1){40}} \put(90,70){\vector(-1,-1){40}}
\put(90,30){\vector(1,0){40}}  \put(130,30){\vector(0,1){40}}   \put(130,70){\vector(-1,-1){40}}  \put(130,70){\vector(0,1){40}} 
 \put(130,110){\vector(-1,-1){40}}  
\put(130,30){\vector(1,0){40}}  \put(170,30){\vector(0,1){40}} \put(170,70){\vector(-1,-1){40}} \put(170,70){\vector(0,1){40}}  
 \put(170,110){\vector(-1,-1){40}} \put(170,110){\vector(0,1){40}}  \put(170,150){\vector(-1,-1){40}}  
\put(90,70){\vector(1,0){40}}  \put(130,70){\vector(1,0){40}} \put(170,70){\vector(1,0){40}} 
\put(130,110){\vector(1,0){40}}  \put(170,110){\vector(1,0){40}}  
\put(170,150){\vector(1,0){40}}
\put(50,30){\circle*{3}}
\put(90,30){\circle*{3}} \put(90,70){\circle*{3}} 
\put(130,30){\circle*{3}}  \put(130,70){\circle*{3}}  \put(130,110){\circle*{3}} 
\put(170,30){\circle*{3}}  \put(170,70){\circle*{3}}  \put(170,110){\circle*{3}}  \put(170,150){\circle*{3}} 
\put(210,30){\circle*{3}}  \put(210,70){\circle*{3}} \put(210,110){\circle*{3}} \put(210,150){\circle*{3}} \put(210,190){\circle*{3}}
  
\put(220,30){\dots}
\put(40,10){(0,0)}
\put(80,10){(1,0)}
\put(120,10){(2,0)}
\put(160,10){(3,0)}
\put(200,10){(4,0)}
\end{picture}
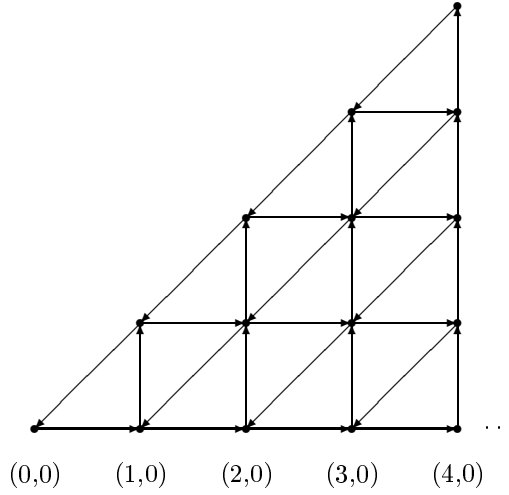
Figure 7. The connected component of the graph for $N=2$.
\end{center}

Suppose we are given a function $\xi \in L^2(\RR^{\ZZ}, d\nu_0^{gc})$ with 
Fourier expansion $\xi =\sum_{I \in \mathcal{I}} \hat{\xi}_I H_I$. 
We can actually turn our directed graph into a weighted graph by assigning to 
each directed edge $(o_1,o_2)$ the Fourier coefficient $\hat{\xi}_I$ corresponding
to the unique multi-index $I$ such that $I \in o_1$ and $\tau I \in o_2$.
Note that each Fourier coefficient will be assigned to one and only one
edge and each edge will have assigned one and only one Fourier coefficient, 
since there is a one-to-one correspondence between the edges of our graph 
and the set $\mathcal{I}$ of multi-indices. For example the edge $(o(I), o(\tau I))$
will have attached the weight $\hat{\xi}_I$.
\begin{center}
\begin{picture}(400,50)(-50,0)
\put(125,40){$\hat{\xi}_I$}
\put(50,30){\circle*{3}}
\put(200,30){\circle*{3}}
\put(50,30){\vector(1,0){150}}
\put(45,10){$o(I)$}
\put(195,10){$o(\tau I)$}
\end{picture}
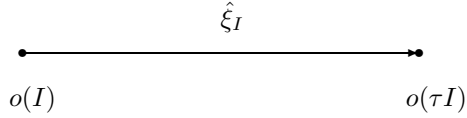
Figure 8. A directed edge and its attached weight.
\end{center}
      
 It is interesting to note that if $\xi$ is $Y_0$-closed then the weights of the
graph discussed before sum up to zero along any directed cycle of the graph except 
the loop of the connected component corresponding to $N=0$. Indeed the closedness 
condition of $\xi$ imposes no restriction on the coefficient $\hat{\xi}_0$.
Also note that the $Y_0$-closedness condition (\ref{feq4})
$$\hat{\xi}_{I+\delta_{(n+1)}}-\hat{\xi}_{I+\delta_n} = \hat{\xi}_{\tau^n(I+\delta_{1})}-
\hat{\xi}_{\tau^n(I+\delta_{0})} \quad n \in \ZZ, \quad I \in \mathcal{I}$$
plus the square-integrability of $\xi$ are equivalent to the property that 
the weights of the graph associated to $\xi$ sum up to zero around any 
directed cycle of any connected component corresponding to $N \geq 1$. 
However if $\xi$ is $Y_0$-exact then $\hat{\xi}_0=0$ and hence the 
weights of the graph sum up to zero around any directed cycle of the graph.
If $\xi$ is $Y_0$-exact and is constructed as in lemma \ref{lem2} then we can say that
in any connected component of the graph all but finitely many weights are zero.

The above can be explained from an algebraic topologic point of view. 
We can turn our directed graph into a $2$-dimensional $\Delta$-complex 
(see A. Hatcher's book on Algebraic Topology, \cite{Hat}) by attaching
enough discs to cycles of the graph such that each of the connected components
$N \geq 1$ can be retracted to a single point. We do not attach a disc 
onto the loop of the connected component $N=0$. After the attaching 
process the $2$-dimensional $\Delta$-complex can be retracted 
to the disjoint union of a circle with a countable number of points. 
The Fourier coefficients of a $Y_0$-exact function form a coboundary 
of our $2$-dimensional $\Delta$-complex and the Fourier coefficients 
of a $Y_0$-closed function form a cocycle for our $2$-dimensional 
$\Delta$-complex. Since the cohomology group $H^1(C, \mathbb{R})$ of 
a circle $C$ is one-dimensional we expect the space of $Y_0$-exact functions
to have codimension one inside the space of $Y_0$-closed functions. 

\section*{Acknowledgments}
 
I would like to thank Professor George A. Elliott for his suggestion to use a Fourier analysis 
approach to solve the problem discussed in this paper.
% and for welcoming me at the Fields Institute for Research in Mathematical Sciences during part of this research.
%I would also like to 
%thank him for his guidance, patience and never ending encouragement. 

\end{document}